\documentclass[a4paper,11pt,reqno]{amsart}
\usepackage{amsmath,amssymb,amsthm}
\usepackage{epsfig}
\usepackage{amsmath,amssymb,amsthm}
\newtheorem{thm}{Theorem}

\newtheorem{lem}[thm]{Lemma}

\theoremstyle{definition}

\newtheorem{rem}{Remark}

\title{An Extension of hibi's palindromic theorem}

\author[D. Lee]{Daeseok Lee}
\address{Daeseok Lee \\ Korea Science Academy of KAIST\\ Pusan 614-822, Korea}
\email{apig5942@naver.com}

\author[H.-K. Ju]{Hyeong-Kwan Ju}
\address{Hyeong-Kwan Ju \\ Department of Mathematics \\ Chonnam National University \\ Gwangju 500-757, Korea}
\email{hkju@jnu.ac.kr}

\date{\today}
\begin{document}
\date{\today}
\begin{abstract}
Hibi showed that the polynomial in the numerator of the Ehrhart series of a reflexive polytope
is palindromic. We proved that those in the numerator of the Ehrhart series of every graph polytope (defined later) of
the bipartite graph is palindromic. From this, one of the conjectures
(raised in the A205497 of OEIS \cite{[O]}) follows immediately.
\end{abstract}
\maketitle

\section{Introduction}

 In the paper \cite{[BJY]}, B\'{o}na and his coauthors enumerated certain weighted graphs with the following conditions :
For a given positive integer $k$, a nonnegative integer $n$ and a simple graph
$G=(VG, EG)$ with $VG=[k]$, where $[k] := \{1, 2, \cdots, k\}$ and $[k]_* :=[k]\cup\{0\}$,
 we consider the set\\
 $$W(n;G) := \{\alpha=(n_1, \cdots, n_k)\in ([n]_*)^k\mid ij\in EG \Rightarrow n_i+n_j
  \leq n\}. $$\\
 We call $\alpha$ satisfying the conditions in the set given above \textit{(vertex-)
 weighted graph}. The question is what the number of weighted graphs would be.
  Of course it depends on the type of graphs. From the question above, we can come up with the following geometric object immediately and naturally.
 Again, let $G = (VG,EG)$  be a simple graph with $VG=[k]$. Then the {\bf graph polytope} $P(G)$ associated with the graph $G$ is defined as follows:\\
$$P(G) := \{(x_1, x_2, \cdots, x_k)\in [0,1]^k \mid ij\in EG \Rightarrow x_i + x_j \leq 1 \}$$\\
Recently, Lee and Ju(\cite{[LJ]}) and Ju, Kim and Seo(\cite{[JKS]}) showed many interesting facts about the volume of the graph polytopes for various types of graphs.
 Main result in this paper is that the numerator polynomial of Ehrhart series of every graph polytope of the finite simple bipartite graph  is
palindromic.
 We introduce a few definitions, necessary lemmas, Hibi's palindromic theorem (or HPT shortly) and main results  in section 2.
 Basically, we will follow the proof scheme of  HPT.
 In the last section we introduce an example (in fact, a conjecture raised in the A205497 of OEIS) where our main theorem applied.
 Also some related remarks were included there.

\section{Main Results}

We will use same notations as those in the book of Beck and Sinai \cite{[BS]}. Let $S$ be an integral convex polytope in ${\mathbb{R}}^d$.
Then we call $L_S(t)=|tS \cap {\mathbb{Z}}^d|$ {\it the Ehrhart polynomial}
of $S.$ From now on, if we mention a polytope, then it is always convex. If the origin belongs to the interior of the polytope $S$, then its {\it dual polytope} $S^*$
is defined as $S^* = \{ x \in {\mathbb{R}}^d \mid  x \cdot y \le 1 \text{ for all }   y \in S \}.$
A polytope $S$ is {\it reflexive} if and only if both $S$ and its dual $S^*$
are integral polytopes. {\it Ehrhart series} of $S$ is the generating function
$Ehr_S(z):=1+ \sum_{t \ge 1} L_S(t)z^t.$ A polynomial $p(z)$ of degree $d$ is  {\it palindromic} if $z^d p(1/z)=p(z).$\\

\begin{thm} {\tt\bf HPT}(\cite{[H]} or \cite{[BS]})
For a given integer $d-$polytope $P$ containing the origin in its interior, if its Ehrhart series is the rational function
of the form
$$Ehr_P(z)= \frac{h_d z^d + h_{d-1}z^{d-1} + \cdots + h_1z + h_0}{(1-z)^{d+1}},$$
then the polytope is reflexive if and only if $h_k=h_{d-k}$ for all $0 \leq k \leq \frac{d}{2}$ (that is, the polynomial in the numerator of the Ehrhart series is palindromic).\\
\end{thm}

Now, let $P$ be a $d-$dimensional integer polytope as follows:
\begin{equation}\label{prp}
P=\{ x=(x_1,x_2,\cdots,x_d)^t  \in {\mathbb{R}}^d :x \ge 0, Mx \le u \},
\end{equation}
where $M$  is the integer matrix and  $u=(1,1,\cdots,1)^t$.
We call the polytope given as the set (\ref{prp}) above a positive reflexive polytope and $M$ a representing matrix of $P$.
If, moreover, all the row sums of the matrix $M$ in (\ref{prp}) are same, we call $P$ a regular positive reflexive polytope.

\begin{lem}
Let $P$ be a regular positive reflexive polytope.Then there exists a positive integer $k \geq 2$ such that
$Q=kP-u$ is a (ordinary) reflexive polytope. Moreover, $P=Q \cap {\mathbb{R}}_{\ge 0}^d$.
\end{lem}

Proof: Let $M=(a_{ij})$ be a representing matrix of size $m \times d$. Let $k=a_{1,1}+a_{1,2}+ \cdots + a_{1,d}+1 \leq 2$.
Then the system of inequalities $x_i \leq 0$ and $a_{i,1}x_1+a_{i,2}x_2+ \cdots a_{i,d}x_d \leq k$
which constructs $kP$ becomes $-(x_i-1) \leq 1$ and  $a_{i,1}(x_1-1)+a_{i,2}(x_2-1)+ \cdots a_{i,d}(x_d-1) \leq 1.$
That is,
\begin{equation}\label{qq}
 Q=kP-u = \{ x \in {\mathbb{R}}^d :
 \begin{bmatrix}
 M\\
 -I_d
 \end{bmatrix}
 x \le u \}.
 \end{equation}

 Facts that this $Q$ in \ref{qq} is a reflexive polytope and
 $P=Q \cap {\mathbb{R}}_{\ge 0}^d$ can be seen easily from the matrix construction. \qed \\

We call will the new polytope obtained as in \ref{qq} the {\it dilated reflexive polytope} of $P$.

\begin{thm}
Let $P$ be a regular positive reflexive polytope and let $k$ be a positive integer obtained from the previous lemma.
Then, for integer $t$,
 \begin{equation*}
 L_{p^o}(t)= \begin{cases}
 0 & \text{if  $t<k$} \\
 L_p(t-k) & \text{if $t \geq k$ } \\
 \end{cases}
 \end{equation*}
\end{thm}

Proof: Let's denote the set of all integer lattice points in a set $X$ by $I[X]$.
Since $I[kP^o]= \{ u \}$ and $u \not\in I((k-1)P^o)$, the case $t<k$ is obvious. The claim is that
$I[u+(t-k)P]=I[tP^o]$. Let $Q$ be the dilated reflexive polytope of the polytope $P$.
Since $u+(t-k)P \subset kP^o + (t-k)P =tP^o$, we have $I[u+(t-k)P] \subset I[tP^o]$. For the other inclusion,
\begin{eqnarray*}
I[tP^o] & \subset & I[(u+(t-k+1)Q^o) \cap (u+ {\mathbb{R}}_{\ge 0}^d)]\\
        & =       & I[(u+(t-k)Q) \cap (u+ {\mathbb{R}}_{\ge 0}^d)] \\
        & =       & I[u+(t-k)P].
\end{eqnarray*}
The inclusion above holds because $tP =kP+(t-k)P \subset u+Q + (t-k)Q = u+(t-k+1)Q$.
The first equality comes from the fact that $I[tQ^o]=I[(t-1)Q]$ (from the proof of HPT) for reflexive polytope $Q$ and integer $t$,
and so does the second one from the fact that $P=Q \cap {\mathbb{R}}_{\ge 0}^d$. \qed \\

Next theorem can be shown by using the Stanley's reciprocity laws. Its proof follows the scheme of the HPT proof.(Refer \cite{[BS]},\cite{[H]},\cite{[RS]} or \cite{[FK]}.)
\begin{thm} \label{pal}
Let $P$ be a $d-$dimensional regular positive reflexive polytope.
Then the $h^* - $polynomial for $P$ is a palindromic polynomial
of degree $d-k+1$.
\end{thm}

Proof:
 \begin{eqnarray*}
 \frac{h_0 z^d+h_1 z^{d-1}+ \cdots +h_d}{(1-z)^{d+1}} & = & (-1)^{d+1} \frac{1}{z} Ehr_P ( \frac{1}{z})\\
                                                      & = & \frac{1}{z} Ehr_{P^o} (z)\\
                                                      & = & \frac{1}{z} \sum_{t=1}^{\infty} L_{P^o}(t) z^t \\
                                                      & = & \frac{1}{z} \sum_{t=k}^{\infty} L_{P}(t-k) z^{t-k}z^k \\
                                                      & = & z^{k-1} Ehr_{P}(z)\\
                                                      & = &  \frac{h_d z^{d+k-1} +h_{d+k-2} z^{d+k-2}+ \cdots +h_0 z^{k-1}}{(1-z)^{d+1}}.
 \end{eqnarray*}
By comparing the coefficients of numerator polynomials in the first and last terms in the previous identity obtained, we can derive that
$h_i=0$ for $i>d-k+1$ and $h_j=h_{d-k+1-j}$ for $0 \leq j \leq d-k+1$.    \qed \\

For a given bipartite graph $G=(VG, EG)$ with $VG=\{ v_1, v_2, \cdots, v_d \}$ and $EG= \{ e_1, e_2, \cdots, e_m \}$,
the {\it incidence matrix} $B=B(G)=(b_{ij})$ of the graph $G$ is the $d \times m$ matrix
whose rows are indexed by $VG$, whose columns are indexed by $EG$, and where
\begin{equation*}b_{ij}= \begin{cases}1 & \text{ if a vertex } v_i \text{ is incident with an edge } e_j,\\
                       0 & \text{ otherwise. }\\
          \end{cases}
\end{equation*}
It is obvious that the graph polytope $P(G)$ can be written as follows:
\begin{equation*}
P(G)=\{ x=(x_1,x_2,\cdots,x_d)^t:x \ge 0, B^tx \le u \},
\text{ where } u=(1,1,\cdots,1)^t \in {\mathbb{R}}^d .
\end{equation*}
Note that all the column sums of the matrix $B$ are 2.
All the above statements (i.e., Lemmas and Theorems) work through for $M:=B^t, P:=P(G), Q:=Q(G),$ and, for example, $k=3$.
(Note here that the graph $G$ should be bipartite.(Refer \cite{[BJY]} for details.)
Thus, next theorem holds in an obvious manner. The proof is almost mimic of that of Theorem  \ref{pal}.

\begin{thm} \label{pal1}
Let $P=P(G)$ be a graph polytope for a given bipartite graph $G$.
Then the $h^* - $polynomial for $P$ is a palindromic polynomial
of degree $d-2$.
\end{thm}


\section{An Application and Further Remarks}

Let the linear graph $L_n$ be the graph $L_n =(VL_n, EL_n)$, where $VL_n=[n]$ and $EL_n=\{ i(i+1) \mid i \in [n-1]\}.$
Then its graph polytope $P(L_n)$ is given as in the following:
$$P(L_n) := \{(x_1, x_2, \cdots, x_n)\in [0,1]^n \mid  x_i + x_{i+1} \leq 1(i=1,2,\cdots,n-1) \}$$
For the graph polytope $P_n=P(L_n),$ set $L_{P_n}(t):=|tP_n \cap {\mathbb{Z}}^n |$ for $t=1,2,\cdots.$ The Ehrhart series for this $P_n$ is
$$Ehr_{P_n} (z) :=1+ \sum_{t \ge 1} L_{P_n} (t)z^t . $$

In the paragraph describing a sequence of an id number A205497 in OEIS(\cite{[O]}) five Conjectures were listed, namely, Conjecture1 through Conjecture5.
Among them, Conjecture1 is of our concern. Anti-diagonal sequences in A205497 given below

$$\begin{array}{ccccccc}
 1,     & 1,     &    1,  &   <1>, &   1,   &      1, & \cdots\\
1,      & 3,     &  <7>,  &  14,   &    26, &     46, & \cdots\\
1,      &<7>,    &  31  , &  109,  &   334, &    937, & \cdots\\
<1>,    &   14,  & 109,   & 623,   &  2951, &  12331, & \cdots\\
1,      &26,     & 334,   & 2951,  & 20641, & 123216, & \cdots\\
1,      &46,     &937,    &  12331,& 123216,& 1019051,& \cdots\\
\vdots  & \vdots & \vdots & \vdots & \vdots & \vdots  & \vdots
 \end{array}$$

\noindent are the coefficients of numerator polynomials of the rational functions of the form
$Ehr_{P_n} (z)=\frac{h_n (z)}{(1-z)^{n+1}} (n=0,1,2,\cdots)$ given by the ordinary generating functions,
for example, $$\frac{1}{(1-z)^3}, \frac{1+z}{(1-z)^4}, \frac{1+3z+z^2}{(1-z)^5}, \frac{<1>+<7>z+<7>z^2+<1>z^3}{(1-z)^6}, \cdots,$$
of Ehrhart polynomial  $L_{P_n} (t)$ in $t$ given row sequences of A050446 in OEIS which is written as below:

$$\begin{array}{ccccccc}
1,      & 1,     &  1,   &   1,   &   1,   &    1,  & \cdots\\
1,      & 2,     &  3,   &   4,   &   5,   &    6,  & \cdots\\
1,      & 3,     &  6,   &  10,   &  15,   &    21, & \cdots\\
1,      & 5,     & 14,   &  30,   &  55,   &    91, & \cdots\\
1,      & 8,     & 31,   &  85,   & 190,   &   371, & \cdots\\
(1),    & (13),  & (70), &  (246),& (671), &  (1547), & \cdots\\
\vdots  & \vdots & \vdots & \vdots & \vdots & \vdots  & \vdots
 \end{array}$$

For example, $$\begin{array}{lll}
Ehr_{P_5} (z)&= & \frac{h_5 (z)}{(1-z)^{6}}= \frac{1+7z+7z^2+z^3}{(1-z)^6}\\
             &= &1+13z+70z^2+246z^3+671z^4+1547z^5+\cdots.\end{array}$$

Conjecture1 states that the polynomials $h_n (z)$ are palindromic for all $n.$ However, the conclusion follows easily
by the Theorem \ref{pal}  since the polytope $P_n$ is a graph polytope
for the simple graphs $L_n$ for all nonnegative integer $n$.

\begin{rem}
It is well-known that $h_n(1)=E_n$, the volume of graph polytope
is $vol(P_n)=\frac{E_n}{n!}$, where $E_n$ is the $n-$th Euler number(cf.A000111 in OEIS), and that $\sum_{n \ge 0} vol(P_n)z^n=secz+tanz$.
(See Elkies\cite{[E]} or Ju et. al.\cite{[JKS]}.) We can generate $E_n$ by the recursion formula
$E_{n+1}=\frac{1}{2} \sum_{k=0}^n \binom{n}{k} E_k E_{n-k}, \: E_0=0.$
\end{rem}

\begin{rem}
Generating function $k_m(y)$ of the $m$th column sequence of A050446 in OEIS can be written as follows(See \cite{[BJY]} for more details.):\\

$$\begin{array}{ll}
k_0(y)&=[-y,1]   \\
k_1(y)&=[-y,y,1]= \cfrac{1+y}{1-y-y^2}  \text{ (generating function of the Fibonacci sequence)}\\
k_2(y)&=[-y,y,-y,1] \\
k_3(y)&=[-y,y,-y,y,1] \\
k_4(y)&=[-y,y,-y,y,-y,1] \\
      & \vdots
\end{array}$$
where
$$[a_0, a_1, a_2, \cdots, a_r]= \cfrac{1}{a_0 + \cfrac{1}{a_1 +\cfrac{1}{a_2 + \cfrac{1}{\ddots+\cfrac{1}{a_r}}}}}.$$
\end{rem}

It is obvious that
\begin{equation*}
\sum_{i \geq 0} k_i(y)z^i = G(y,z) = \sum_{j \geq 1} Ehr_{P_n}(z)y^j.
\end{equation*}
However, the shape of the bivariate generating function $G(y,z)$ is not known.

\begin{rem}
Let $C_t=(c_{ij})$ be a square $t \times t$ matrix defined as follows:
$$c_{ij}= \begin{cases}1 & \text{ if } i+j \le t+1,\\
                       0 & \text{ otherwise. }\\
          \end{cases}  $$
For a given matrix $M$, we also let $s(M):= \text{the sum of all entries of the matrix of }M.$
Then $L_{P_n}(t)=s((C_t)^n).$ (Refer \cite{[BJY]} for details.)
\end{rem}

\begin{rem}
We can define the graph polytope associated with hypergraph $H=(V,E)$, the corresponding Ehrhart polynomial and Ehrhart series.
In a similar way, we can show that main theorems still hold. It might be interesting to apply our methodology
to the graph polytopes associated with $k-$uniform hypergraphs.
\end{rem}


\end{document}